\documentclass{article}

\usepackage{amsmath,amssymb}
\usepackage[all,knot]{xy}
\usepackage{color}

\def\mapr#1{\smash{\mathop{\buildrel{#1}\over\longrightarrow}}}

\def\mapd#1{\Big\downarrow\rlap{$\vcenter{\hbox{$#1$}}$}}

\newtheorem{theorem}{Theorem}
\newtheorem{definition}{Definition}

\def\proof{{\bf Proof.}}
\def\qed{\hfill{$\square$}}
\def\cA{{\cal A}}

\def\cC{{\cal C}}
\def\cS{{\cal S}}
\def\rg{{\frak g}}
\def\id{{\hbox{\bf Id}}}

\title{Generalization of the Sullivan construction for Transitive Lie Algebroids}
\author{Mishchenko, Alexander S. (Moscow, Russia)\thanks{Partly supported by RFBR 11-01-00057-a, 10-01-92601-KO\_a, 11-01-90413-Ukr\_f\_a and State program RNP 2.1.1.5055} and \\
Ribeiro, Jose (Braga, Portugal)}

\begin{document}

\maketitle
\section{Preliminaries and Formulation of the Problem}

D.Sullivan \cite{S-77} (see also the book by H.Whitney \cite{whit-1957}) considered a new model for underlying cochain complex for classical cohomologies with rational coefficients for arbitrary simplicial spaces that gives the isomorphism with classical rational cohomologies (p.297, Theorem 7.1).
We
apply the key ideas
developed by K.MacKenzie \cite{Mck-05} and J.Kubarski \cite{Kub-91} to a generalization of the D.Sullivan model for transitive Lie algebroids. The generalization is based on the existence of the inverse image of the transitive Lie algebroids and on the property of transitive Lie algebroids being trivial over contractible manifolds. Using these properties one can construct an underlying cochain complex of differential forms on simplicial space.

Namely, for transitive Lie algebroid ${\cal A}\rightarrow TM\rightarrow M$ one can define an underlying Sullivan model as graded algebra
$\Omega^{*}_{ps}({\cal A};M)$ of piecewise smooth differential forms
on smooth manifold $M$ with smooth combinatorial structure.

Each piecewise smooth differential form
$\omega\in\Omega^{p}_{ps}({\cal A};M)$ is a collection \linebreak
$\omega=\{\omega_{\sigma}\in \Omega^{p}(\cal A;\sigma)\}$ where $\sigma$ runs the collection of all closed simplices and if $\phi:\tau\subset\sigma$ is a face then
$
\phi^{!}(\omega_{\sigma})=\omega_{\tau}
$.

One has the natural morphism
$$
\Omega^{*}({\cal A};M)\rightarrow \Omega^{*}_{ps}({\cal A};M)
$$
of differential algebras.
In the case of smooth combinatorial manifold we prove that the cohomologies of the Sullivan model for cochain complex are isomorphic to usual cohomologies with coefficients of the Lie algebroid.

\subsection{Definition of transitive Lie Algebroid}

Let $\cA$ be a transitive  Lie algebroid. These means that we have a
smooth manifold $M$, tangent bundle $TM$
$$
    \begin{array}{c}
    TM\\ \mapd{p}\\M
    \end{array}
$$
a vector bundle $\cA$ over $TM$,
$$
    \begin{array}{c}
    \cA\\\mapd{\gamma}\\TM\\ \mapd{p}\\M
    \end{array}
$$
or
$$
    \xymatrix{
      \cA\ar[r]_{\gamma}\ar[d]_{p_{A}}&TM\ar[d]_{p}\\
      M\ar[r]_{=}&M
    }
$$
and an operation $[\cdot,\cdot]$ on the space of sections
$\Gamma^{\infty}(\cA;M)$ (or $\cS_{\hbox{ec}}^{\infty}(\cA;M)$ ) that satisfies the following conditions:
\begin{enumerate}
\item The anchor $\gamma$ is fiberwise surjective;
\item The operation $[\cdot,\cdot]$ forms a Lie algebra structure on $\Gamma^{\infty}(\cA;M)$;
\item The anchor $\gamma$ induces the Lie algebra homomorphism
$$
    \gamma_{\Gamma}:\Gamma^{\infty}(\cA;M)\mapr{}\Gamma^{\infty}(TM;M);
$$
\item Action of the algebra $\cC^{\infty}(M)$ on $\Gamma^{\infty}(\cA;M)$ satisfies the natural condition:
$$
    [\xi,f\cdot\eta]=\gamma_{\Gamma}(\xi)(f)\cdot \eta + f\cdot[\xi,\eta].
$$
\end{enumerate}

A differential form $\Phi$ of degree $j$ is a skew-symmetric polylinear map of
$\cC^{\infty}(M)$--modules
$$
    \Phi:\underbrace{\Gamma^{\infty}(\cA;M)\times\Gamma^{\infty}(\cA;M)
    \times\cdots\times\Gamma^{\infty}(\cA;M)}_{\hbox{j times}}
    \mapr{}\cC^{\infty}(M)
$$

We denote by $\Omega^{j}(\cA;M)$ the family of all differential forms of degree $j$.

\subsection{Restriction to Submanifold}

Let $\varphi:N\hookrightarrow M$ be a compact submanifold (possibly with boundary).
We can define the restriction of the Lie algebroid $\cA$ over $TM$ on the submanifold $N$:

$$
    \xymatrix{\varphi^{!!}(\cA)=&
      (T\varphi)^{*}(\cA)\ar[r]\ar[d]_{\gamma_{TN}}&  \cA\ar[d]_{\gamma}\\
      &TN\ar[r]^{T\varphi}\ar[d]_{p_{N}}&TM\ar[d]_{p_{M}}\\
      &N \ar[r]^{\varphi}&M
    }
$$

It is evident that, in the case in which $\varphi:N\hookrightarrow M$ is an open subset,  the restriction $\varphi^{!!}(\cA)$
has a natural structure of the transitive Lie algebroid.

For arbitrary submanifold $\varphi:N\hookrightarrow M$, the structure  of the Lie algebroids on the $\varphi^{!!}(\cA)$ is defined by natural extension of sections from $\Gamma^{\infty}(\varphi^{!!}(\cA);M)$ to $\Gamma^{\infty}(A;M)$.

There is natural homomorphism of modules:
$$
    \Omega^{j}(\cA;M)\mapr{\varphi^{\cA}}\Omega^{j}(\varphi^{!!}(\cA);N)
$$

The homomorphism $\varphi^{\cA}$ commutes with the exterior differentials on differential forms.

Remark. When $\varphi:N\hookrightarrow M$ is an open subset then without lost of correctness we say that
$\Gamma^{\infty}(\varphi^{!!}(\cA),N)\approx\Gamma^{\infty}(\cA,N)$ and
$\Omega^{j}(\varphi^{!!}(\cA);N)\approx\Omega^{j}(\cA;N)$

\subsection{The Mayer-Vietoris sequence for
forms on Lie algebroids}

Let $\cA$ be a Lie algebroid on a smooth
$M$. Assume that $U$ and $V$ are open subsets of $M$ such that $M=U\cup
V$. Put $W=U\cap V$. One has the short exact sequence

\begin{equation*}
\begin{array}{ccccccccc}
0 & \mapr{} & \Omega^{*}(\cA;M) & \mapr{} &
\Omega^{*}(\cA;U)\times \Omega^{*}(\cA;V) &
\mapr{} & \Omega^{*}(\cA;W) & \mapr{} & 0
\end{array}
\end{equation*}
which gives the long exact sequence in homologies:

\begin{equation*}
\begin{array}{l}
\cdots \mapr{} H^{p-1}(\cA;W)
\mapr{}  H^{p}(\cA;M)\mapr{}
H^{p}(\cA;U)\times
H^{p}(\cA;V)\mapr{}\\\\
\mapr{}  H^{p}(\cA;W)
\mapr{}  H^{p+1}(\cA;M)\mapr{}H^{p+1}(\cA;U)\times
H^{p+1}(\cA;V)  \mapr{}
\cdots
\end{array}
\end{equation*}

\vspace{5mm}

This long sequence is the Mayer-Vietoris sequence in Lie
algebroids.

\section{Piecewise smooth forms}
Let $M$ be a smooth manifold with a smooth fixed simplicial combinatorial structure.
\begin{definition}
A piecewise differential form $\omega\in\Omega^{p}_{ps}(\cA;M)$ is a collection of the differential forms
$\{\omega_{\sigma}\in\Omega^{p}(\cA; \sigma)\}$ such that for any face $\phi_{\tau,\sigma}:\tau\hookrightarrow\sigma$ one has
$$
\phi_{\tau,\sigma}^{\cA}(\omega_{\sigma})=\omega_{\tau},
$$
where
$$
\phi^{\cA}_{\tau,\sigma}:\Omega^{p}(\cA;\sigma)\mapr{}\Omega^{p}(\cA;\tau)
$$
is the homomorphism generated by inclusion $\phi_{\tau,\sigma}$.
\end{definition}

The space $\Omega^{*}_{ps}(\cA;M)=\bigoplus\limits_{p=0}^{\infty}
\Omega^{p}_{ps}(\cA;M)$ is the graded differential algebra with differential $$
d_{ps}:\Omega^{p}_{ps}(\cA;M)\mapr{}\Omega^{p+1}_{ps}(\cA;M).
$$
The homology is denote by
$$
H^{*}_{ps}(M;\cA)= H(\Omega^{*}_{ps}(\cA;M),d_{ps}).
$$
There is the natural map
$$
r^{\cA}:\Omega^{p}(\cA;M)\mapr{}\Omega^{p}_{ps}(\cA;M)
$$
defined by
$$
r^{\cA}(\omega)=\{r^{\cA}(\omega)_{\sigma}\}, \quad r^{\cA}(\omega)_{\sigma}=\phi_{\sigma}^{\cA}(\omega),
$$
where $\phi_{\sigma}:\sigma\hookrightarrow M$ is inclusion of the simplex
$\sigma$.
It is evident that $r^{\cA}$ commutes with $d$:
$$
r^{\cA}\circ d=d_{ps}\circ r^{\cA}.
$$
So one has the homomorphism
$$
H(r^{\cA}):H^{*}(M;\cA)\mapr{}H^{*}_{ps}(M).
$$

\begin{theorem}
The map $H(r^{\cA}):H^{*}(M;\cA)\mapr{}H^{*}_{ps}(M;\cA)$
is the isomorphism.
\end{theorem}

\proof
The proof consist of two steps.
The first step is to extend the exact Mayer-Vietoris sequence to piecewise
smooth differential forms. Assume that open set $U\subset M$ is the union
of some open simplices. Then one can extend the notion of the piecewise smooth differential forms $\omega\in\Omega^{p}_{ps}(\cA;U)$ as a collection
$\{\omega_{\sigma}\in\Omega^{p}(\cA; \sigma_{U}): \sigma_{U}=\sigma\cap U \}$
which satisfies the property: if $\tau\subset\sigma$ is the face,
$\phi_{\tau,\sigma}:\tau_{U}\hookrightarrow\sigma_{U}$ is natural inclusion, then

$$
\phi_{\tau,\sigma}^{\cA}(\omega_{\sigma})=\omega_{\tau}.
$$

For example, the star of the (open) simplex $\sigma$, $\star(\sigma)$, is the open subset that is union of open simplices. So we have the commutative diagram
\xymatrix{
\cdots\ar[r]& H^{p-1}(U\cap V;\cA)
\ar[r]\ar[d]^{H(r^{\cA})}&
H^{p}(U\cup V;\cA)\ar[r]\ar[d]^{H(r^{\cA})}&
{\begin{array}{c}
H^{p}(U;\cA)\\\oplus\\H^{p}(V;\cA)
\end{array}}
\ar[r]\ar[d]^{H(r^{\cA})}&
\cdots\\
\cdots\ar[r]& H^{p-1}_{ps}(U\cap V;\cA)
\ar[r]&
H^{p}_{ps}(U\cup V;\cA)\ar[r]&
{\begin{array}{c}
H^{p}_{ps}(U;\cA)\\\oplus\\H^{p}_{ps}(V;\cA)
\end{array}}
\ar[r]&
\cdots
}

Using the five lemma one can reduce the proof to the $U=\star(\sigma)$.

The second step is to prove the theorem for special case $U=\star(\sigma)$.
The manifold $U$ is diffeomorphic to the disc and therefore the transitive
Lie algebroid $\cA$ on $U$ is trivial (\cite{Mck-05}, theorem 7.3.18, page 286). Hence
$$
\Omega^{*}(\cA;U)\approx\Omega^{*}(U)\otimes\Omega^{*}(\rg)
$$
with differential $d^{\cA}=d^{TM}\otimes\id+\id\otimes d^{\rg}$.
For piecewise smooth differential forms one has similar properties:
$$
\Omega^{*}_{ps}(\cA;U)\approx\Omega^{*}_{ps}(U)\otimes\Omega^{*}(\rg)
$$
with differential $d^{\cA}_{ps}=d^{TM}_{ps}\otimes\id+\id\otimes d^{\rg}$.
Then on homology level one has
$$
H^{*}(U;\cA)\approx H^{*}(U;R)\otimes H^{*}(\rg),
$$
$$
H^{*}_{ps}(U;\cA)\approx H^{*}_{ps}(U;R)\otimes H^{*}(\rg),
$$
Hence the map $H(r^{\cA})$ for $U$
$$
H^{*}(U;\cA)\approx H^{*}(U;R)\otimes H^{*}(\rg)\mapr{H(r^{\cA})}
H^{*}_{ps}(U;R)\otimes H^{*}(\rg)\approx H^{*}_{ps}(U;\cA)
$$
has the form
$$
H(r^{\cA})=H(r^{TM})\otimes\id
$$
So the proof is reduced to the Sullivan theorem (\cite{S-77}).
\vspace{10mm}\qed


\begin{thebibliography}{aaaa}
\bibitem{whit-1957}
Whitney, H., \emph{Geometric Integration Theory}, Princeton
University Press, 1957.
\bibitem{S-77}Dennis Sullivan, \textit{Infinitesimal computations in topology.} Publ.
I.H.E.S. 47 (1977) 269-331.

\bibitem{Mck-05}
Mackenzie, K.C.H., \textit{General Theory of Lie groupoids and Lie algebroids}, Cambridge University Press, Cambridge, 2005.

\bibitem{Kub-91}Kubarski, J.,
  \textit{The Chern-Weil homomorphism of regular Lie algebroids},
  Publications du Department de Mathematiques, Universite Claude Bernard - Lyon-1,
  (1991) {4--63}.
\end{thebibliography}
\end{document}